# MIXED NON-EXPANSIVE AND POTENTIALLY EXPANSIVE PROPERTIES OF A CLASS OF SELF-MAPS IN METRIC SPACES

M. De la Sen

Institute of Research and Development of Processes. Campus of Leioa, Bizkaia, SPAIN

email: manuel.delasen@ehu.es

**Abstract**- This paper investigates self-maps $T: X \to X$ which satisfy a distance constraint in a metric space which mixed point-dependent non-expansive properties, or in particular contractive ones, and potentially expansive properties related to some distance threshold. The above mentioned constraint is feasible in certain real-world problems.

## 1. Introduction

Fixed point theory and related techniques are of increasing interest for solving a wide class of mathematical problems where convergence of a trajectory or sequence to some equilibrium set is essential. Recently, the subsequent set of more sophisticated related problems are under strong research activity:

1) In the, so-called, p-cyclic non-expansive or contractive self-maps map each element of a subset $A_i$ of an either metric or Banach space **B** to an element of the next subset $A_{i+1}$ in a strictly ordered chain of p subsets of **B** such that $A_{p+1} = A_1$. If the subsets do not intersect then fixed points do not exist and their potential relevance in Analysis is played by best proximity points, [1-2]. Best proximity points are also of interest in hyperconvex metric spaces, [3-4].

2) The so-called Kannan maps are also being intensively investigated in the last years as well as their relationships with contractive maps. See, for instance, [5-6], [11].

3) Although there is an increasing number of theorems about fixed points in Banach or metric spaces, new related recent results have been proven. Some of those novel results are, for instance, the generalization in [7] of Edelstein´s fixed point theorem for metric spaces by proving a new theorem. Also, an iterative algorithm for searching a fixed point in nonexpansive mappings in Hilbert spaces has been proposed in [8]. On the other hand, an estimation of the size of an attraction ball to a fixed point has been provided in [9] for nonlinear differentiable maps.

4) Fixed point theory can be also used successfully to find oscillations of solutions of differential or difference equations which can be themselves characterized as fixed points. See, for instance, [9-10], [12-13].

This manuscript is devoted to investigate self-maps $T: X \to X$ in a metric space $(X, d)$ which satisfy the constraint $d(Tx, Ty) - d(x, y) \leq -K d(x, y) + M$; for some real constants $K \geq 0$, $M \geq 0$. It is direct to see that $d(Tx, Ty) \leq d(x, y)$; i.e. $T: X \to X$ is non-expansive, if $d(x, y) \geq M/K$. Also,

$$d(x, y) < M/K \Rightarrow d(Tx, Ty) \leq (1-K) d(x, y) + M < M/K; \quad \forall x, y \in X \qquad (1.1)$$

Then, the self-map $T:X \to X$ exhibits the following constraint under (1.1) provided that it is continuous: $T:A_{xy} \to A_{\omega z}$ where $A_{xy} \subset X$ is the open circle of center $c_{xy} \in X$ of radius $R:=M/K$ for each given $x,y \in A_{xy}$ and $A_{\omega z} \subset X$ is an open circle of center at some $c_{\omega z} \in X$ also of radius R. Note that $A_{xy}$ can be distinct from $A_{\omega z}$. However, if $T:X \to X$ is not continuous then the existence of the above circles is not ensured but only that (1.1) holds. Note that (1.1) does not guarantee that, contrarily to the case of large distances fulfilling $d(x,y) \geq M/K$ if $K \neq 0$ and $d(x,y) = \infty$ if $K=0$, the self-map $T:X \to X$ cannot be guaranteed to be non-expansive, while it can be eventually expansive, for small distances fulfilling $d(x,y) < M/K$ ; $x,y \in X$. The objective of this paper is the investigation of self-maps $T:X \to X$ which such mixed properties related to some distance threshold.

## 2. Basic Distance Property and Related Motivating Example

Let $(X,d)$ be a metric space and T a self-map from X to X. Such a self-map is uncertain in the sense that the distance is subject to the following constraint:

$$d(Tx,Ty) - d(x,y) \leq -Kd(x,y) + M \; ; \; \forall x,y \in X \text{, some real constants } K \geq 0, M \geq 0 \quad (2.1)$$

In order to discuss the feasibility of (2.1), note the following:

1) If $M=0$ and $K \in (0,1]$ then (2.1) is the usual contractive constraint of Banach contraction principle and $T:X \to X$ is strictly contractive. If $K=M=0$ then $T:X \to X$ is non-expansive. If $M=0$, $K=1$ and the inequality in (2.1) is strict for $x,y (\neq x) \in X$ then $T:X \to X$ is weakly contractive.

2) If K=1 then $d(Tx,Ty) \leq M$ ; $\forall x,y \in X$. Since T is a self-map on X, the validity of the constraint (2.1) is limited to the set family $\hat{A}_T := \{A_i \subset X : (diam(A_i) \leq M \wedge T(A_i) \subset A_j ; some A_j \in \hat{A}_T)\}$ of bounded subsets of X. In this case, $d(T^j x, T^j y) \leq M$ ; $\forall j \in \mathbb{Z}_+$ provided that $x,y \in A_\alpha \in \hat{A}_T$ and T maps X to some member $A_i$ of $\hat{A}_T$ for each given $x,y \in X$. In other words, the image of T is restricted as $T:X \to X | A_i$ (for some $A_i \in \hat{A}_T$ which depends, in general, on x and y ) so that $d(Tx,Ty) \leq M$ in order to (2.1) to be feasible, i.e. $Tx,Ty$ are in some set of the family $\hat{A}_T$ if the pair x, y in X is such that $d(x,y) > M$. Note that $T:X \to X$ is not necessarily a retraction from X to some element of $\hat{A}_T$ since $T(A_i) \subseteq A_j$ for $A_i, A_j(\neq A_i) \in \hat{A}_T$. Note that $T:X \to X | A_i$ can possess a fixed point if K=1 and ((2.1) holds.

3) If $K>1$ then $d(Tx,Ty) \leq M$ if $x=y$ ; $x,y \in X$, and

$$d(x,y) \geq M/(K-1) \Rightarrow 0 \leq d(Tx,Ty) \leq d(x,y) - \frac{M}{K-1} < d(x,y) \text{ if } x,y(\neq x) \in X$$

Then if $x, y \in X$ exist such that $d(x, y) \in \left(0, \dfrac{M}{K-1}\right)$ then (2.1) is impossible for any self-map T on X since it would imply $d(Tx, Ty) < 0$. For $x = y$, (2.1) holds for self-maps T on X such that $d(Tx, Ty) \leq M$. Fixed points can exist only in trivial cases as, for instance, $X := \left\{x : d(x, y) \geq \dfrac{M}{K-1}; \forall y \in X\right\}$ is a set of isolated points with a minimum pair-wise distance threshold so that $T : X \to X$ is such that $T(y) = x \in X$; $\forall y \in X$.

4) The case of interest discussed through this paper for (2.1) is when $M > 0$ and $K \in [0, 1)$. It is shown that the self-map $T : X \to X$ exhibits contractive properties for sufficiently large distances which exceed a minimum real threshold while it might possibly be expansive for distances under such a threshold. A related motivating example follows.

*2.1* *Example 2.1*: Note that (2.1) is equivalent to:

$$d(Tx, Ty) \leq (1-K) d(x, y) + M \; ; \; \forall x, y \in X \text{, for some } M > 0 \qquad (2.2)$$

Eq. 2.1 is relevant, for instance, in the following important physical problem. Let a linear time-invariant n-th order dynamic system be:

$$\dot{x}(t) = A x(t) + \eta_x(t) \qquad (2.3)$$

with $A \in \mathbf{R}^{n \times n}$ being a stability matrix whose fundamental matrix satisfies $\|e^{At}\| \leq K_0 e^{-\alpha_0 t}$; $\forall t \geq 0$ for some positive real constants $K_0$ (being norm-dependent) and $\alpha_0$ and $\eta : [0, \infty) \times X \to \mathbf{R}^n$ being an unknown uniformly bounded perturbation of essential supremum bound satisfying $\mathrm{ess} \sup\limits_{\infty > t \geq 0} \|\eta_x(t)\| \leq M_0 < \infty$; $\forall x \in X$. The unique solution of (2.3) for $x(0) = x_0$ is:

$$x(t) = e^{At} x_0 + \int_0^t e^{A(t-\tau)} \eta_x(\tau) d\tau \qquad (2.4)$$

Direct calculation with (2.4) for the norm-induced distance $d(x, y) := \|x - y\|$; $\forall x, y \in X$ yields:

$$d(x(t), y(t)) = \|x(t) - y(t)\| \leq K_0 e^{-\alpha_0 t} \|x_0 - y_0\| + \dfrac{K_0}{\alpha_0} \sup_{0 \leq \tau < \infty} \|\eta_x(\tau) - \eta_y(\tau)\|$$

$$\leq (1 - K) d(x_0, y_0) + M \; ; \; \forall t \geq h_0 := \dfrac{1}{\alpha_0} \ln \dfrac{K_0}{1 - K} \qquad (2.5)$$

with $\infty > M \geq \dfrac{2 K_0 M_0}{\alpha_0}$, $K := 1 - K_0 e^{-\alpha_0 h_0} \in (0, 1)$. Now, let $X \subset \mathbf{R}^n$ the state space of (2.1), generated by (2.4), subject to $x_0 \in X$ and $(X, d)$ is a complete metric space. Define the state transformation $T_h x(kh) = x[(k+1)h]$ on X which generates the sequence of states $\{x(kh)\}_{k=0}^{\infty}$ being in X if $x_0 \in X$ with h being any real constant which satisfies $h \geq h_0$. Then, the self-map $T_h : X \to X$ satisfies (2.1). Note that the system (2.3) is always globally Lyapunov stable for any bounded initial conditions in view of (2.5). If the perturbation is identically zero then the origin is globally asymptotically Lyapunov stable since A is a stability matrix. This follows also from (2.5) since the self-map $T_h$ on X is a contraction which has zero as its unique fixed and equilibrium point so that

$x(kh+\tau) = e^{A\tau} x(kh) \to 0$ as $k \to \infty$; $\forall \tau \in [0,h)$; $\forall h \geq h_0$. Thus, $x(t) \to 0$ as $t \to \infty$.

However, in the presence of the perturbation, the origin is not globally asymptotically stable (although the system is globally stable) and it exhibits ultimate boundedness since for sufficiently large distances $d(T_h x(kh), T_h y(kh)) \geq \dfrac{M}{K}$ (respectively, $d(x(kh), y(kh)) > \dfrac{M}{K}$), the self-map is non-expansive (respectively, contractive). Then, $0 \leq d(T_h x(kh), T_h y(kh)) \leq d(x(kh), y(kh))$ respectively, $d(T_h x(kh), T_h y(kh)) < d(x(kh), y(kh))$ ). But such properties are not guaranteed if $d(x(kh), y(kh)) < \dfrac{M}{K}$ which can lead to $T_h : X \to X$ being expansive. □

*Remark 2.2*: Example 2.1 emphasizes the fact that some real-world problems exist where certain self-maps T from X to X are neither contractive nor expansive everywhere in X while such a map is guaranteed to be contractive for sufficiently large distances between any two points in X exceeding a known real threshold. For small distances, the self-map could be potentially expansive, or, as in the dynamic system of Example 2.1, unclassified as expansive, non-expansive or contractive. In Example 2.1, this last situation is due to the presence of unknown perturbations of known prescribed upper-bound. Note that in Example 2.1, the self-map from X to X is guaranteed to be point-wise contractive or potentially expansive for each given pair in X accordingly to the distance between them. Note that the global asymptotic Lyapunov´s stability relies on a contractive mapping and a zero equilibrium point which is a global attractor and a fixed point of a certain mapping from initial conditions to subsequent points of the state- trajectory solution. However, this is not a requirement for global Lyapunov´s stability where the above mapping may be relaxed to be non- contractive or to the so-called ultimate boundedness (which also implies global stability with eventual local instability around the equilibrium) related to boundedness of the state trajectory solution for arbitrary bounded initial conditions. In this case, the above mapping can be locally expansive and then contractive for trajectory – solution points being sufficiently far away from the equilibrium point.

### 3. Main Results

This section is devoted to formalize the general context of the described problem to the light of Fixed Point Theory. A first main result follows:

**Theorem 3.1**. Assume that $K \in (0,1)$ and consider any bounded set $X_0 \subset X$ with $diam(X_0) \leq R$

**(i)** Assume that $R \geq M/K$. Then, the restricted map $T|X_0$ of T from $X_0$ to X is non-expansive (i.e. $d(Tx, Ty) \leq d(x, y)$ ) for any pair $x, y \in X_0$ such that $d(x, y) \geq M/K$ and weakly contractive (i.e. $d(Tx, Ty) < d(x, y)$ ) for any pair $x, y \in X_0$ such that $d(x, y) > M/K$.

**(ii)** The distance between the iterates $T^j x$ and $T^j y$ is uniformly bounded ; $\forall x, y \in X_0$, $\forall j \in \mathbf{Z}_+$ and there exists a bounded subset $X_1$ of X fulfilling $X_1 \supset X_0 \neq X_1$ such that $T^j x, T^j y \in X_1$; $\forall j \in \mathbf{Z}_+$.

**(iii)** If Property (ii) holds then all iterate $T^j x$ enters a compact convex subset $X_\alpha$ of $X$ ; $\forall x \in X_0$; $\forall j \geq j_0$ and some finite integer $j_0$; i.e. the sequence $T^j x$ is permanent; $\forall x \in X_0$. Also, any two pairs of iterates $T^j x, T^j y$ enter within a compact subset of $X_\alpha$ of prescribed diameter $\frac{M}{K} + \varepsilon$; $\forall x, y \in X_0$, $\forall j \geq j_0$ and some finite integer $j_0$.

*Proof:* *(i)* It follows by direct inspection from (2.2).

*(ii)* Direct recursive calculation with (2.2) for $j \in \mathbf{Z}_+$ and any bounded $X_0 \subset X$ with $diam(X_0) \leq R$ yields:

$$d(T^j x, T^j y) \leq (1-K)^j d(x, y) + M \sum_{i=0}^{j-1} (1-K)^{j-1-i} \qquad (3.1)$$

$$\leq (1-K)^j d(x, y) + \frac{M}{K}\left(1-(1-K)^j\right) \leq d(x, y) + \frac{M}{K} \leq R + \frac{M}{K} < \infty \ ; \ \forall j \in \mathbf{Z}_+ \qquad (3.2)$$

, $\forall x, y \in X$ since $K \in (0,1)$.

*(iii)* From the first inequality of (3.2),

$$\limsup_{j \to \infty} d(T^j x, T^j y) \leq \frac{M}{K}\left(1-(1-K)^j\right) \leq \frac{M}{K} < \infty \ ; \ \forall x, y \in X \qquad (3.3)$$

so that for any given real constant $\varepsilon > 0$, it exists a positive integer $j_{01} = j_{01}(\varepsilon)$ such that from (3.2)

$$d(T^j x, T^j y) \leq \frac{M}{K} + (1-K)^{j_{01}}\left(R - \frac{M}{K}\right) \leq \frac{M}{K} + \varepsilon \ ; \ \forall j \geq j_{01} \ ; \ \forall x, y \in X \qquad (3.4)$$

provided that $\varepsilon \geq (1-K)^{j_{01}}\left(R - \frac{M}{K}\right)$, or equivalently, $j_{01} \geq \ln\frac{\varepsilon}{R - M/K} - \ln|1-K|$ provided that $\frac{M}{K} \leq R \leq \frac{\varepsilon}{|1-K|} + \frac{M}{K}$. Now, assume that $R > \frac{\varepsilon}{|1-K|} + \frac{M}{K}$. Then, from (3.3) and the definition of limit superior, there exists a finite positive integer $j_{02} = j_{02}(\delta, R)$ for any arbitrary given positive real constant $\delta$ such that $d(T^j x, T^j y) \leq \frac{M}{K} + \delta$ ; $\forall j \geq j_{02}'$. Then, choose $\delta = \frac{\varepsilon}{|1-K|}$ for the given $\varepsilon$. Thus, (3.4) holds ; $\forall j \geq j_{02} := j_{01} + j_{02}'$. Finally, assume that $0 \leq R < \frac{M}{K}$. Then, from (3.3), it exists $j_{03} = j_{03}(\varepsilon)$ such that $d(T^j x, T^j y) \leq \frac{M}{K} + \varepsilon$ ; $\forall j \geq j_{03}$. As a result, for any bounded set $X_0 \subset X$ with $diam(X_0) \leq R$, it exists a finite positive integer $j_0 = j_0(\varepsilon, R)$ such that $d(T^j x, T^j y) \leq \frac{M}{K} + \varepsilon$; $\forall x, y \in X_0$, $\forall j \geq j_0$, for some finite integer $j_0$. Thus, for each given real $\varepsilon > 0$, there is a compact convex subset $X_\alpha \supset X_0$ of X where all the iterates $T^j x$ enter; $\forall j \geq j_0$, for some finite integer $j_0$. Furthermore, any two iterates $T^j x, T^j y$ are within a compact convex subset of $X_\alpha$ of prescribed diameter $\frac{M}{K} + \varepsilon$; $\forall x, y \in X_0$, $\forall j \geq j_0$. □

*Remark 3.2.* Note that Theorem 3.1 (i) does not conclude that the self-map $T: X_0 \to X_0$ is expansive for some pair $x, y \in X_0$ (i.e. $d(Tx, Ty) > d(x, y)$) if $d(x, y) < M/K$ but only that the upper-bound $(1-K)d(x, y) + M$ of $d(Tx, Ty)$ is upper-bounded by $d(x, y)$. Thus, $d(x, y) < M/K$ for some pair $x, y \in X_0$ is a necessary (but not sufficient) condition for $T: X_0 \to X_0$ to be expansive for that pair.

□

*Remark 3.3.* Note that $X_0$ in Theorem 3.1 is not required to be convex. Theorem 3.1 (ii) guarantees that $T^j x, T^j y \in X_1 \supset X_0$ although eventually it may not belong to $X_0$. □

It is now of interest to characterize in some sense a subset $X_e$ of X such that the restricted map $T|X_e$ is a self-map from $X_e$ to $X_e$ which satisfies the constraints:

$$K_1 d(x, y) \leq d(Tx, Ty) \leq \min\left((1-K)d(x, y) + M, K_2 d(x, y)\right); \forall x, y \in X_e \subseteq X \quad (3.5)$$

for some real constants $K \in [0, 1)$ $M > 0$, $K_1 > \max(1-K, 0)$, $K_2 \geq K_1$. This will allow later on the definition of subsets of X where the self-map T is contractive, expansive or non-expansive. It would be proven later on (see Corollary 3.5) that (3.5) is impossible everywhere in X if $K_1 > 1$.

**Theorem 3.4.** Assume that $K \in [0, 1)$. Then, there is a family of nonempty bounded subsets of X for which (3.5) holds and, then trivially, a subfamily of nonempty bounded convex subsets of X with the same property.

*Proof*: The constraints (3.5) are guaranteed under two possibilities for each $x, y \in X_e(x) \subseteq X$ where $X_e(x) := \{y \in X : (3.5) \text{ holds}\}$ is a point-dependent subset of X, namely:

$$K_1 d(x, y) \leq d(Tx, Ty) \leq (1-K)d(x, y) + M \leq K_2 d(x, y); \forall x, y \in X_e(x) \subseteq X, \text{ some } M > 0 \quad (3.6)$$

which implies:

$$K_1 d(x, y) \leq d(Tx, Ty) \leq K_2 d(x, y) \leq (1-K)d(x, y) + M; \forall x, y \in X_e(x) \subseteq X \quad (3.7)$$

The constraint (3.6) is subject to the necessary conditions:

$$d(x, y) \in \left[\frac{M}{K + K_2 - 1}, \frac{M}{K + K_1 - 1}\right]; \forall x, y \in X_e(x) \quad (3.8)$$

$$\frac{M K_1}{K + K_2 - 1} \leq K_1 d(x, y) \leq d(Tx, Ty)$$

$$\leq (1-K)d(x, y) + M \leq K_2 d(x, y) \leq \frac{M K_1}{K + K_1 - 1}; \forall x, y \in X_e(x)$$

$$\Rightarrow d(Tx, Ty) \in \left[\frac{MK_1}{K + K_2 - 1}, \frac{MK_1}{K + K_1 - 1}\right]; \forall x, y \in X_e(x) \quad (3.9)$$

Since T is a self-map on $X$, any pair $x, y \in X_e(x)$ has to satisfy simultaneously (3.8)-(3.9) so that

$$d(x, y) \in \left[ \frac{M}{K + K_2 - 1}, \frac{M}{K + K_1 - 1} \right] min(1, K_1) ; \forall x, y \in X_e(x) \quad (3.10)$$

under the constraint (3.6). Since $K_2 \geq K_1$, the constraint (3.7) requires

$$d(x, y) \leq \frac{M}{K + K_2 - 1} ; \forall x, y \in X_e$$

$$d(Tx, Ty) \leq K_2 d(x, y) \leq \frac{MK_2}{K + K_2 - 1} \leq \frac{MK_1}{K + K_1 - 1} ; \forall x, y \in X_e \quad (3.11)$$

The last inequality of (3.11) follows directly if $K_1 > 1 - K$ since $K_2 \geq K_1 > 1 - K$ implies that

$$\frac{K + K_1 - 1}{K_1} = 1 - \frac{1 - K}{K_1} \geq 1 - \frac{1 - K}{K_2} = \frac{K + K_2 - 1}{K_2} \quad (3.12)$$

Combining (3.11)-(3.12), one gets that (3.7) holds if

$$d(x, y) \in \left[ 0, \frac{M}{K + K_2 - 1} \right] min(1, K_2) ; \forall x, y \in X_e(x) \quad (3.13)$$

Thus, it is clear the existence of a countable family of nonempty bounded subsets $\{X_{ei}(x)\}$ of $X_e(x); \forall x \in X$ defined by

$$X_{ei}(x) := \left\{ y \in X : d(x, y) \leq \frac{M}{K + K_1 - 1} min(1, K_1) \right\} \subset X_e(x) \subset X ; \forall x \in X \quad (3.14)$$

since $X_e(x) := \{y \in X : (3.5) \ holds\} \supset X_{ei}(x); \forall x \in X$. From the above developments, it turns out that there exists a convex subset in the family $\{X_{ei}(x)\}$ which is convex and then a subfamily of the set $\{X_{ei}\}$ which possess such a property. □

Theorems 3.1 and 3.4 lead to the following important conclusion:

**Corollary 3.5**. Assume that $K \in (0,1)$. Then, the following properties hold if (3.5) holds:

(i) If $max(1-K, 0) < K_1 \leq K_2 \leq 1$ then $T : X \to X$ is non-expansive.

(ii) If $max(1-K, 0) < K_1 \leq K_2 < 1$ then $T : X \to X$ is (strictly) contractive and then it has a fixed point.

(iii) If $K_1 \in [0,1)$ and $K_2 > 1$ then the restriction of T to $\hat{X}(x)$, $T | \hat{X}(x) := (T : X | \hat{X}(x) \to X)$, $\forall x \in X$, is non-expansive where $\hat{X}(x) := \left\{ y \in X : d(x, y) \geq \frac{M}{K} \right\} \subset X ; \forall x \in X$ but $T | \hat{X}_{ei}(x)$ is weakly contractive for all sets $X_{ei}(x)$ defined in (3.14) resulting to be

$$X_{ei}(x) := \left\{ y \in X : d(x, y) \leq \frac{M}{K + K_1 - 1} \right\}$$ since $K_1 > 1 ; \forall x \in X$. As a result $T : X \to X$ is neither contractive nor expansive on X.

(iv) If $K_2 \geq K_1 = 1$ then $T : X \to X$ is not contractive, and

$$K_1 d(x,y) \leq d(Tx,Ty) \leq (1-K)d(x,y)+M \leq max\left(K_2 d(x,y), \frac{K_1 M}{K+K_1-1}\right); \forall x,y \in X$$

(3.15)

If $K_1 > 1$ then neither (3.6) nor (3.7) is feasible for any $x, y \in X$ and (3.5) is not feasible either; $\forall x, y \in X$.

*Proof*: Properties *(i) –(ii)* follow rom Theorem 3.1. Property (iii) follows from Theorem 3.4 since:

$K_1 > 1 \Rightarrow M/K \geq M/(K+K_1-1) \Rightarrow \hat{X}(x) \cap X_{ei}(x) = \varnothing ; \forall x \in X$

Then, for any $x, y \in X$, if $y \in \hat{X}(x)$ then $y \notin \hat{X}_{ei}(x)$ and conversely. The constraints (3.15) follow directly from (3.5) and its necessary condition $d(x,y) \leq \frac{M}{K+K_1-1}$ ; $\forall x, y \in X$. It is now proven by contradiction that neither (3.6) nor (3.7) is feasible for all given pair x, y in X if $K_1 > 1$. A necessary condition for (3.6) to hold for each $x, y \in X$ is that $d(x,y) \in \left[\frac{M}{K+K_2-1}, \frac{M}{K+K_1-1}\right]$. Thus, $T: X \to X$ is not expansive which contradicts $d(x,y) < K_1 d(x,y) \leq d(Tx,Ty)$ ; $\forall x, y(\neq x) \in X$ if $K_1 > 1$. Also, if (3.7) holds; $\forall x, y \in X$ then $d(x,y) \leq \frac{M}{K+K_2-1}$ which leads to the same above contradiction if $K_1 > 1$. On the other hand, a necessary condition for (3.5) to hold is that

$K_1 d(x,y) \leq (1-K)d(x,y) + M \Rightarrow d(x,y) \leq \frac{M}{K+K_1-1} < \infty ; \forall x, y \in X$

which contradicts $K_1 > 1$. *Property (iii)* has been proven. *Property (iv)* follows from (3.6) and *Property (iii)* for $K_2 \geq K_1 = 1$. □

## ACKNOWLEDGMENTS

The author thanks to the Spanish Ministry of Education by its support of this work through Grant DPI2009-07197. He is also grateful to the Basque Government by its support through Grant GIC07143-IT-269-07.

## REFERENCES

[1] S. Karpagam and S. Agrawal, "Best proximity point theorems for p-cyclic Meir- Keeler contractions", *Fixed Point Theory and Applications*, Vol. 2009, Article ID 197308, 9 pages, doi:10.1155/FPTA/2009/197308.

[2] A.A. Eldred and P. Veeramani, "Existence and convergence of best proximity points", *Fixed Point Theory*, Vol. 323, No. 2, pp. 1001-1006, 2006.

[3] A. Amini-Harandi, A.P. Farajzadeh, D. O. Reagan and R.P. Agarwal, "Coincidence point best approximation and best proximity theorems for condensing set-valued maps in hyperconvex metric spaces", *Fixed Point Theory and Applications*, Vol. 2009, Article ID 543154, doi:10.1155/FPTA/2009/543154.

[4] A. Amini-Harandi, A.P. Farajzadeh, D. O. Reagan and R.P. Agarwal, "Best proximity pairs for upper semicontinuoyus set-valued maps in hyperconvex metric spaces ", *Fixed Point Theory and Applications*, Vol. 2009, Article ID 648985, doi:10.1155/FPTA/2009/648985.

[5] M. Kikkawa and T. Suzuki, "Some similarity between contractions and Kannan mappings", *Fixed Point Theory and Applications*, Vol. 2008, Article ID 649749, 8 pages, doi:10.1155/FPTA/2008/649749.

[6] Y. Enjouji, N. Nakanishi and T. Suzuki, "A generalization of Kannan´s fixed point theorem", *Fixed Point Theory and Applications*, Vol. 2009, Article ID 192872, 10 pages, doi:10.1155/FPTA/2009/192872.

[7] T. Suzuki, "A new type of fixed point theorem in metric spaces", *Nonlinear Analysis: Theory, Methods & Applications,* Vol. 7, No. 11, pp. 5313-5317, 2009.

[8] Q YU Liu, W Y Zeng and NJ Huang, "An iterative method for generalized equilibrium problems, fixed point problem and variational inequality problems", *Fixed Point Theory and Applications*, Article ID 531308, 20 pages, doi:1155/2009/531308, 2009.

[9] E. Catinas, "Estimating the radius of an attraction ball ", *Applied Mathematics Letters*, Vol. 22, No. 5, pp. 712-714, 2009.

[10] M.R. Pournaki and A. Razani, "On the existence of periodic solutions for a class of generalized forced Lienard equations", *Applied Mathematics Letters*, Vol. 20, No. 3, pp. 248-254, 2009.

[11] M. De la Sen, "Some combined relations between contractive mappings, Kannan mappings, reasonable expansive mappings, and T- stability", *Fixed Point Theory and Applications*, Vol. 2009, Article ID 815637, 25 pages, doi:10.1155/JAM/2009/815637.

[12] M. De la Sen, "Some fixed point properties of self-maps constructed by switched sets of primary self-maps on normed linear spaces", *Fixed Point Theory and Applications,* vol. 2010, Article ID 438614, 20 pages, doi: 10.1155/2010/438614.

[13] M. De la Sen, "About robust stability of dynamic systems with time delays through fixed point theory", *Fixed Point Theory and Applications,* vol. 2008, Article ID 480187, 20 pages, doi: 10.1155/2008/480187.